\documentclass[11pt]{article}
\usepackage[text={6.75in,9.5in},centering,letterpaper]{geometry}
\usepackage{amsfonts}
\usepackage{amsmath}
\usepackage{amsthm}
\usepackage{amssymb}
\usepackage{amscd}
\usepackage{enumerate}
\usepackage{hyperref}
\usepackage{algorithm}
\usepackage{algorithmic}
\usepackage{graphicx}

\newcommand{\beq}{\begin{equation}}

\newcommand{\pd}[2]{\dfrac{\partial #1}{\partial #2}}

\newcommand{\eeq}{\end{equation}}
\newcommand{\ut}{\tilde{u}}

\def\k0{\kappa_0}

\newtheorem{thm}{Theorem}[section]

\theoremstyle{definition}
 
\theoremstyle{Remark}

\numberwithin{equation}{section}

\title{Global Error Analysis and Inertial Manifold Reduction\thanks{The research was support in part by NSF grants DMS-1115408 and DMS - 1419047 and
University of Kansas Undergraduate Research Awards.}}
\author{Yu-Min Chung\thanks{Department of Mathematics, College of William and Mary, Williamsburg, VA 23187 {\tt (ychung@wm.edu)}.}\and
        Andrew Steyer\thanks{Department of Mathematics, University of Kansas, Lawrence, KS 66045 {\tt (a518s941@ku.edu)}.}\and
        Michael Tubbs\thanks{Department of Mathematics, University of Kansas, Lawrence, KS 66045 {\tt (tubbs.irishcraze@gmail.com)}.}\and
        Erik S. Van Vleck\thanks{Department of Mathematics, University of Kansas, Lawrence, KS 66045 {\tt (erikvv@ku.edu)}.}\and
        Mihir Vedantam\thanks{Department of Mathematics, University of Kansas, Lawrence, KS 66045 {\tt (vedantam.mihir@gmail.com)}.}
}
        
\begin{document}

\maketitle

\begin{abstract}
Four types of global error for initial value problems are considered in a common framework.
They include classical forward error analysis and shadowing error analysis together with extensions of 
both to rescaling of time. To determine the amplification of the local error that bounds the global error
we present a linear analysis similar in spirit to condition number estimation for linear systems of equations.
We combine these
ideas with techniques for dimension reduction of differential equations via a boundary value formulation of 
numerical inertial manifold reduction.  These global error concepts are exercised to illustrate their utility on
the Lorenz equations and inertial manifold reductions of the Kuramoto-Sivashinsky equation.
\end{abstract}

\vskip -1.0in

\section{Introduction}

In many complex systems modeled using differential equations the slow dynamics drive the system. There is a vast literature on
inertial manifold techniques to determine the mapping between the slow dynamics and the fast dynamics.
This decouples the system and focuses attention on the (often) low dimensional slow dynamics that drive the system. Once such a low dimensional reduction is achieved,
then one would like to infer the behavior of the system from simulations of the reduced equations. An often overlooked problem is 
in assessing whether the global error on the inertial manifold can be controlled and in what sense. The standard approach to global
error analysis is classical forward error analysis for initial value differential equations in which the initial condition is the 
same for both the exact solution and the numerical approximation. Shadowing error analysis generalizes this in a significant way by
allowing for slightly different initial conditions for the exact and approximate solutions. This expands the class of problems for
which long time error statements are possible, from contractive problems to those with a splitting between expansive and contractive
modes. A further refinement that has been investigated in the shadowing literature involves the rescaling of time when differential equations have a non-trivial attractor.

Our contribution in this paper is to develop a unified approach to global error analysis for initial value
problems that can be used to determine when there is uncertainty in the numerical approximation of solutions
of differential equations;
we also show that the technique is applicable in the context of inertial manifold dimension reduction.
We report on initial numerical experiments for a new inertial manifold reduction technique combined with
an assessment of the global error in approximating the reduced equations. The inertial manifold technique which we outline here and describe in
more detail in \cite{CSVV14} involves first performing a time dependent linear decoupling tranformation and then determining the mapping between the
slow and fast dynamics implicitly by solving a boundary value problem. Information obtained during the solution of the boundary value problem
is then employed to assess the relationship between the local error and if possible the global error as characterized by forward or shadowing
error analysis and with or without rescaling of time. In this paper we highlight four different types of long time global error analysis.
We will also exercise recently developed techniques for dimension reduction in differential equations in the context of these types of global error analysis.

Shadowing based techniques for global error analysis involve relaxing the requirement that
the initial conditions for the exact solution and the numerical solution agree. This has the
effect that the linearized error equation need not be solved forward in time. This allows for
positive global error statements for a larger class of problems over long time intervals, i.e., for
problems that are not contractive such as in the case of a system with positive Lyapunov exponents.
Shadowing also provides a framework to allow for rescaling of time,
i.e., allowing for perturbations in the time step,
(see the work of Coomes, Kocak, and Palmer \cite{CKPZAMP}, \cite{CKPNM}, and Van Vleck \cite{VVshad95}, \cite{VVshad00}).
Rescaling of time is especially important when there is a periodic orbit or more general non-trivial attractor \cite{HLW10}.
Work on numerical shadowing includes the ground breaking work of Hammel, Yorke, and
Grebogi \cite{HYG1}, \cite{HYG2}, the work of
Chow, Lin, and Palmer (e.g. \cite{CLP}, \cite{CP2}), the numerical work of Sauer and Yorke \cite{SY},
and the initial work on breakdown of shadowing of Dawson, Grebogi, Sauer, and Yorke \cite{DGSY}.

Inertial manifolds, first introduced by Foias, Sell, and Temam \cite{Foias1988309} for dissipative dynamical systems, are finite dimensional, exponentially attracting, positively invariant Lipschitz manifolds.  Similar concepts are slow manifolds in slow-fast system introduced in meteorology and widely used in weather forecasting \cite{Daley1981450, Tribbia1982, lorenz1986existence, camassa1995geometry}, and center-unstable manifolds in the classic sense.  In fact, \cite{Debussche199173} shows that a slow manifold is a special type of inertial manifold, and as mentioned in the original work \cite{Foias1988309}, it can be described as a global center-unstable
manifold.  The main application is the inertial manifolds reduction, meaning that dynamical systems restricted on the manifolds, whose long-term dynamics coincide with those of the original system without introducing errors.  In particular, since the manifold is finite dimensional, the reduced system is also finite dimensional, whereas the original system may arise from an infinite dimensional system.  Because of its importance, there has been tremendous work in regard to its theory and computation, see e.g. \cite{FNST88, CFNT89, NFT89, MPS88, Jolly89, CLT06, HamedGuoTiti2015} and \cite{FNST88, CFNT89, NFT89, MPS88, Jolly89, CLT06, ChungJolly, HamedGuoTiti2015}, respectively.  Recently, the theory of inertial manifolds has been generalized to non-autonomous dynamical systems \cite{ComputingInertial4, CL97, GV97,koksch2002inertial}, and recently, to random dynamical systems \cite{computing_inertial_sde}, and \cite{chekroun2015approximation} (and the references therein).

We take the approach here of decoupling the time dependent linear part of the equation
using techniques that have proven useful in the approximation of Lyapunov exponents.
We first employ an orthogonal change of variables $Q(t)$ that brings the time dependent
coefficient matrix for the linear part of the equation to upper triangular. Subsequently,
we will compute a change of variables that decouples the linear part.
This then gives us equations of the form considered by Aulbach and Wanner
in \cite{AulbachWanner}.
A similar change of
variables has been employed to justify that Lyapunov exponents and Sacker-Sell spectrum
may be obtained from the diagonal of an upper triangular coefficient matrix
(see section 5 of \cite{DVV02} and sections 4 and 5 of \cite{DVV07}).
The references \cite{DVVEncyc} and \cite{DJVV11} (see also the references therein)
provide a summary and overview
of recent work on approximation of Lyapunov exponents and in obtaining the
orthgonal change of variables $Q(t)$. In this paper we consider the smooth orthogonal change of variables
based upon continuous Householder reflectors as developed in \cite{DV1,DV2}.

This paper is outlined as follows.  We first present a framework for global error analysis in section \ref{Errback} .  Techniques to be employed for non-autonomous inertial manifold reduction are in section \ref{IMback} .  In section \ref{Imp} we outline of methods to determine the amplification of the 
local error that determines the global error. This is followed by details of our dimension reduction implementation based upon time dependent linear decoupling transformation and subsequent solution of the inertial manifold equations 
using a boundary value differential equation solver. Section \ref{Results} contains the results of the technique applied to the three dimensional Lorenz 1963 model and to an inertial manifold reduction of a Galerkin approximation of the Kuramoto-Shivashinky equation.


\section{Framework for Global Error Analysis}\label{Errback}

In this section we present a framework for global error analysis of initial value differential equations.
We will focus our attention on four specific characterizations of global error analysis. The differences
among the characterizations is determined by which variables are allowed to differ between the numerically
computed solution and an exact solution. 
This follows the framework developed for shadowing based error analysis in \cite{VVshad95}.

To make these ideas concrete consider a smooth initial value ODE of
the form
\begin{equation}\label{nonlinODE}
\dot u = f(u,t),\,\, u(t_0) = u_0.
\end{equation} 
If we let $\varphi(u_n,h_n;t_n)$ denote the solution operator that advances the state variable $u_n$, $h_n$ time units from $t_n$,
then the exact solution satisfies (for $t_{n+1} = t_n + h_n$),
\[
u(t_{n+1}) = \varphi(u(t_n),h_n;t_n).
\]

A general approach to global error analysis can be obtained using the setup employed in numerical shadowing.
Outlined below are four measures of the computational error in approximating the solution to an initial value
differential equations. Subsequently, we will apply these ideas to the reduction obtained on the inertial manifold to 
assess to the computational error in approximating solutions to these reduced set of equations.


The idea behind shadowing type global error analysis is to use a numerical approximation of the solution as an initial guess
for a functional Newton-type iteration and show that this converges to a nearby exact solution. 
If we let $x = \{x_n\}_{n=0}^N$ and $h=\{h_n\}_{n=0}^{N-1}$ and define
\[
(G(x,h))_n = x_{n+1} - \varphi(x_n,h_n;t_n),\,\, n=0,1,...,N-1,
\] 
then the linear theory for global error analysis in our framework can be reduced to obtaining bounds on a right inverse of an appropriate derivative of $G$.
We consider four possibilities based upon different sets of variables from the $x_n$ and $h_n$
\begin{enumerate}
\item {\bf Case 1 (Variables are $\{x_n\}_{n=1}^N$)}: This is a standard forward error analyis and requires that the exact and approximate solution
have the same initial conditions.
\item {\bf Case 2 (Variables are $\{x_n\}_{n=0}^N$)}:  This is a standard shadowing error analyis and allows that the exact and approximate solution
have different initial conditions.
\item {\bf Case 3 (Variables are $\{x_n\}_{n=1}^N$ and $\{h_n\}_{n=0}^{N-1}$)}:  This is a forward error analyis with rescaling of time.
\item {\bf Case 4 (Variables are $\{x_n\}_{n=0}^N$ and $\{h_n\}_{n=0}^{N-1}$)}:  This is a shadowing error analyis with rescaling of time.

\end{enumerate}

If we linearize with respect to the variables in each of these four cases, then
we obtain for $X_n := \frac{\partial}{\partial x_n} \varphi(x_n,h_n;t_n)$, $f_n := \frac{\partial}{\partial h_n} \varphi(x_n,h_n;t_n)$,
and $\theta>0$ a scaling factor:
\begin{enumerate}
\item {\bf Case 1:}
$(L \Delta x)_n = \Delta x_{n+1} - X_n \Delta x_n$, $n=1,...,N-1$ and $(L \Delta x)_0 = \Delta x_1$.
\item {\bf Case 2:}
Shadowing Error Analysis: $\{x_n\}_{n=0}^N$ are variables. Then 
$(L \Delta x)_n = \Delta x_{n+1} - X_n \Delta x_n$, $n=0,...,N-1$.
\item {\bf Case 3:}
$(L(\Delta x, \Delta h))_n = \Delta x_{n+1} - X_n \Delta x_n -\theta f_n\Delta h_n$, $n=1,...,N-1$ and 
$(L(\Delta x, \Delta h))_0 = \Delta x_1 - \theta f_0\Delta h_0$.
\item {\bf Case 4:}
$(L(\Delta x,\Delta h))_n = \Delta x_{n+1} - X_n \Delta x_n - \theta f_n\Delta h_n$, $n=0,...,N-1$.

\end{enumerate}

Note that Cases 3 and 4 simplify to Cases 1 and 2, respectively, when $\theta=0$.
In Case 1 $L$ is a square, invertible matrix, while in the other cases
if $L$ is full rank, then $L$ has a right or pseudo inverse of the form $L^\dagger = L^T (LL^T)^{-1}$. 
Using the nonlinear shadowing type global error analysis, a uniform bound on the global error $\epsilon$ 
as an amplificiation of the local error $\delta$ 
may be obtained via 
the following fixed point result, essentially proving convergence of a Newton type
method with a frozen Jacobian to find a zero of $G$ near the numerically computed approximate solution.

\begin{thm}\label{fixedpoint}
Suppose that ${\cal X}, {\cal Y}$ are Banach spaces,
$G: {\cal X}\to {\cal Y}$ is $C^1$, and that there exists a
positive constant $\epsilon_0$, a point $z\in{\cal X}$, and
a linear operator $L:{\cal X}\to {\cal Y}$ with right inverse $L^\dagger$ such that
\begin{equation}\label{ineq1}
\|L^{\dagger}(DG(w)-L)\|\leq \frac{1}{2} {\rm \quad for \quad} \|w-z\|\leq\epsilon_0.
\end{equation}
If, for $0<\epsilon\leq\epsilon_0$,
\begin{equation}\label{ineq2}
\|L^{\dagger}G(z)\|\leq \frac{1}{2}\epsilon,
\end{equation}
then the equation $G(w) = 0$ has a solution with
\begin{equation}
\|w - z\| \leq \epsilon.
\end{equation}
\end{thm}

If this theorem holds, then we have a global error analysis with global error given by
$\epsilon := 2\|L^{\dagger}G(z)\|\leq 2 \|L^{\dagger}\|\cdot \|G(z)\|$. If we employ the
infinity norm in sequence space, then $\epsilon \leq 2\|L^{\dagger}\|_\infty \delta$
where $\delta$ is a uniform bound on the local error. Thus, $2\|L^{\dagger}\|_\infty$
represents the amplification of the local error that gives a bound on the global error.
In particular, the change in the state variables $\|\Delta x_n\| \leq  2\|L^{\dagger}\|_\infty \delta$
and in Cases 3 and 4 the change in the time step $|\Delta h_n| \leq 2\theta \|L^{\dagger}\|_\infty \delta$.

We consider four forms for $L$ and hence $L^{\dagger}$ that all have the potential to 
provide long time global error statements depending upon the dynamics of the problem being
considered. For example, classical forward error results can be obtained for contractive problems (all Lyapunov exponents negative)
such as spatially discretized diffusion equations. Classical shadowing results are obtained for problems with
a mixture of positive and negative Lyapunov exponents but no zero Lyapunov exponents. An example where positive
shadowing results have been obtained is for the forced damped pendulum equation. Forward error analysis with rescaling
of time provides long-time global error statements for problems with a single zero Lyapunov exponent
with the remaining Lyapunov exponents negative. An example of such as problem occurs when approximating a stable periodic orbit.
Shadowing with rescaling of time can provide long-time global error statements when approximating nonlinear autonomous 
differential equations with a bounded solution and a single zero Lyapunov exponent. A specific example is the classical Lorenz
1963 model.

\section{Non-Autonomous Inertial Manifold Reduction}\label{IMback}

The dimension of the matrix $L$ is of the order of the number of time steps times the number of dependent variables.
We consider inertial manifold techniques
to reduce the dimension of this matrix. 
We follow the framework in \cite{AulbachWanner}. Consider the nonautonomous dynamical system
\begin{equation}
\label{equ:original NDS}
\begin{cases}
\dot{x} = A(t) x + F(t, x, y),\quad x(t_0) = x_0,\\
\dot{y} = B(t) y + G(t, x, y), \quad y(t_0) = y_0,
\end{cases}
\end{equation}
where $x$ and $y$ are elements of some Banach spaces $X$ and $Y$, respectively, and $A:\mathbb{R} \rightarrow \mathcal{L}(X)$, $B:\mathbb{R} \rightarrow \mathcal{L}(Y)$, $F:\mathbb{R}\times X \times Y \rightarrow X$, and $G:\mathbb{R}\times X \times Y \rightarrow Y$ are mappings satisfying the following assumptions.

\begin{enumerate}[(H1)]
	\item \label{H1} The mappings $A$ and $B$ are locally integrable and there exists $K\geq 1$ and $\alpha < \beta$ such that the evolution operators $\Lambda_A$ and $\Lambda_B$ of the homogeneous linear equations $\dot{x} = A(t)x$ and $\dot{y} = B(t)y$, respectively, satisfy the estimates
	\begin{align*}
	&\| \Lambda_A(t,s) \| \leq Ke^{\alpha (t-s)} \quad \text{ for all } t\geq s,\\
	&\| \Lambda_B(t,s) \| \leq Ke^{\beta (t-s)} \quad \text{ for all } t\leq s.
	\end{align*}
	\item \label{H2} $F(t,0,0) = 0$ and $G(t,0,0) = 0$ for all $t\in \mathbb{R}$. $H:=(F,G)$ is Lipschitz functions---
	\begin{equation*}
	\| H(t, z_1) - H(t, z_2) \| \leq M \| z_1 - z_2 \|,\quad \forall z_1,\;z_2\in Z.
	\end{equation*}
	\item \label{H3} We need the relation between the linear and nonlinear terms, and it is known for the gap condition:
	\begin{equation}
	\label{equ:classic gap}
	{\beta - \alpha} > {2KM}.
	\end{equation}
\end{enumerate}
Under (H1) and (H2), the result in \cite{AulbachWanner} shows that there exists $\Phi: \mathbb{R}\times X \rightarrow Y$ and $\Psi: \mathbb{R}\times Y \rightarrow X$ such that
\begin{equation}
\label{equ:decoupled NDS}
\begin{cases}
\dot{x} = A(t) x + F(t, x, \Phi(t,x))\\
\dot{y} = B(t) y + G(t, \Psi(t,y), y)
\end{cases}
.
\end{equation}
Under this framework, the family of inertial manifolds is $\mathcal{N}_{t} = \text{graph}(\Psi(t, y))$ and the inertial manifold reduction is
\begin{equation}
	\label{equ:IM reduction}
	\dot{y} = B(t) y + G(t, \Psi(t,y), y).
\end{equation}
Proof and further rigorous properties can be found in \cite{AulbachWanner} and \cite{koksch2002inertial}.  In this article, since we are interested in the computational aspect, 
we sketch the way that we use to find $\Psi$, and for more details, we refer readers to our follow-up work \cite{CSVV14} and \cite{CVV14}.  Recall from \cite{AulbachWanner} that 
\begin{equation}
\label{eq:T map}
\mathcal{T}(\psi, y_0, t_0)(t) = \Lambda_B(t,t_0)y_0 + \int_{t}^{t_0} \Lambda_B(t,s)G(\psi(s))\;ds
-\int_{-\infty}^t \Lambda_A(t,s)\;F(\psi(s))\;ds, \quad \forall t \leq t_0,
\end{equation}
and it can be shown that for given $y_0\in Y$ and $t_0\in \mathbb{R}$, $\mathcal{T}$ has a fixed point in a proper Banach space, denoted by $\psi$, which then defines $\Psi$.  Therefore, finding $\psi$ is essential in this computation.  Moreover, from \eqref{eq:T map} one can show that $\psi(t) =: (x(t),\; y(t))$ is the unique solution of the following boundary value problem (BVP)
\begin{equation}
\label{equ:T map ivp}
\begin{cases}
\dot{x} = A(t)x + F(t,x,y),\quad x(-\infty)= 0, \\
\dot{y} = B(t)y + G(t,x,y), \quad y(t_0) = y_0.
\end{cases}
\end{equation}
The main advantage of \eqref{equ:T map ivp} over \eqref{eq:T map} numerically is that the existed BVP packages can be used.

In the numerical approximation of inertial manifolds the predominant approach has been to start from a linearly decoupled equation.
This is often accomplished using an eigen-decomposition of the linear operator in the problem. The approach we will adopt follows that taken in \cite{CSVV14}, which involves the use of a possibly time dependent solution that we linearize about and then decoupling
the time dependent linear operator obtained from the linearization. In particular, if we write
\[
\dot u = f(t,u) = Df(t,u)u + (f(t,u)-Df(t,u)u) \equiv L(t)u + N(t,u)
\]
where $Df(t,u)$ denotes the derivative of $f(t,u)$ with respect to the $u$ variables, then we decouple the linear operator $L(t)$ using techniques employed in finding stability spectra such as Lyapunov exponents and Sacker-Sell spectrum 
(see \cite{CSVV14} for complete details). 


Start with a given ODE initial value problem
\begin{equation}\label{ivporig}
\left\{\begin{array}{c}
\dot{v}=f(v,t)\\
v(0)=u_0
\end{array}
\right.
\end{equation}
which has solution $v(t;u_0)$.  We can decompose $f(u,t)$ as $f(u,t) = Df(v(t;u_0),t)u + N(u,t) \equiv L(t)u+N(u,t)$.  The initial value problem \eqref{ivporig} can then be expressed as
\begin{equation}\label{ivplinvar}
\left\{\begin{array}{c}
\dot{u}=L(t)u+N(u,t)\\
u(0)=u_0
\end{array}
\right.
\end{equation}
Take a fundamental matrix solution $X(t)$ of $\dot{u} = L(t)u$.  We can factor $X(t)$ as $X(t) = Q(t)R(t)$ where $Q(t)$ is orthogonal and
\begin{equation}\label{Rmat}
R(t) = \left[\begin{array}{cc}A(t) & C(t) \\ 0 & B(t) \end{array} \right]
\end{equation}
 with $A(t) \in \mathbb{R}^{p\times p}$ upper triangular and $B(t) \in \mathbb{R}^{q\times q}$ is full.  Under the changes of variables $u(t) = Q(t)z(t)$ where $z(t) = (x(t),y(t))$ with $x(t) \in \mathbb{R}^p$ and $y(t) \in \mathbb{R}^q$ satisfies
\begin{equation}\label{imform}
\left\{\begin{array}{c}
\dot{x} = A(t)x + F(t,x,y)\\
\dot{y} = B(t)y+G(t,x,y)\\
z(0)=Q(0)^T u_0 = (x_0,y_0)
\end{array}
\right.
\end{equation}
The inertial manifold of the system \eqref{imform} consists of solutions of the boundary value problem 
\begin{equation}\label{bvpform}
\left\{\begin{array}{c}
\dot{x} = A(t)x + F(t,x,y)\\
\dot{y} = B(t)y+G(t,x,y)\\
x(-\infty) = 0, y(0) = y_0
\end{array}
\right.
\end{equation}
See \cite{CSVV14} for more details.  The boundary value problem \eqref{bvpform} forms the basis for our dimension reduction techniques.

\subsection{Differentials of Manifolds}
 
The differential $D\Psi$ can be found as the fixed point of the following operator as shown in \cite{DJRV}
\begin{equation}
\label{eq:DT map}
\mathcal{T}_1(\Delta, y_0, t_0)(t) = \Lambda_B(t,t_0) + \int_{t}^{t_0} \Lambda_B(t,s)DG(\psi)\cdot\Delta \;ds
-\int_{-\infty}^t \Lambda_A(t,s)\;DF(\psi)\cdot\Delta\;ds, \quad \forall t \leq t_0,
\end{equation}
and similar to \eqref{bvpform} for \eqref{eq:T map}, we can show that the fixed point of \eqref{eq:DT map} satisfies the following BVP
\begin{equation}
\label{equ:bvp for differential unstable}
\begin{cases}
\dot{\Delta_1} = A(t)\Delta_1 + DF(\psi)\cdot\Delta,\quad \Delta_1(-\infty) = 0\\
\dot{\Delta_2} = B(t)\Delta_2 + DG(\psi)\cdot\Delta,\quad \Delta_2(0) = I_Y		
\end{cases}
,
\end{equation}
where $\psi$ is the solution of \eqref{bvpform}, and $\Delta(t) = (\Delta_1(t), \Delta_2(t))$ is a linear operator from $Y$ to $Z$.  The differential of unstable manifolds is $D\Psi = \Delta_1(0)$.

To implement it, one could couple the \eqref{equ:bvp for differential unstable} with \eqref{bvpform}:
\begin{equation}
\begin{cases}
\dot{x} = A(t)x + F(t,x, y), \quad x(-\infty) = 0 \\
\dot{y}	= B(t)y + G(t,x, y), \quad y(0) = y_0 \\
\dot{\Delta_1} = A(t)\Delta_1 + DF(t,x,y)\cdot\Delta,\quad \Delta_1(-\infty) = 0\\
\dot{\Delta_2} = B(t)\Delta_2 + DG(t,x,y)\cdot\Delta,\quad \Delta_2(0) = I_Y
\end{cases}
.
\end{equation}
Note that $\Delta_1$ is of dimension $\dim(X)\times \dim(Y)$ and $\Delta_2$ is of dimension $\dim(Y)\times\dim(Y)$.

Our focus in this work is on obtaining an estimate of $\|L^\dagger\|_\infty$ in these four
cases when applied to the differential equations obtained by restricting to the inertial manifold.

In our case we will focus on the linear operator generated by the linearization around the solution on the lower dimensional decoupled equation
\[
\dot y = B(t) y + G(t,\Psi(y,t),y) \equiv G_1(t,y)
\]
with flow denoted by $\varphi(t,y)$.  The the linearized equation becomes
\beq\label{LVEonIM}
\dot Y = B(t)Y + G_y(t,\Psi(y,t),y)Y + G_x(t,\Psi(y,t),y)\Psi_y(y,t)Y.
\eeq
Then we can determine the $X_n$ and $f_n$ as $X_n:= X_n(t_{n+1})$ where $X_n(t)$ satisfies
(\ref{LVEonIM}) for $t_n \leq t\leq t_{n+1}$ with $X_n(t_n)=I$ and $f_n \approx G_1(t_{n+1},y(t_{n+1}))$.  We want to avoid solving a large system of equations to form $X_n$.  To do so, we approximate the derivative of $G_1(t,y)$ with respect to $y$ directly with the flow $\varphi$, which we may approximate by time-stepping along the inertial manifold.  We have, for $\Delta x > 0$ and $\Delta t> 0$ small, the approximations
$$\frac{\partial}{\partial y} G_1(t,y)e_j \approx \frac{1}{\Delta x}(G_1(t, y+(\Delta x) e_j)-G_1(t,y)), \quad G_1(t,y)=\dot{y}(t) \approx \frac{1}{\Delta t}(\varphi(\Delta t;y)-\varphi(0;y)), \quad 1,\hdots,p.$$
where $e_j$ is the $j^{th}$ standard basis vector of $\mathbb{R}^p$, where $p=\text{dim}(Y)$.  We make use of this approximation in the next section where we discuss the implementation.

\section{Implementation}\label{Imp}

In this section we describe the techniques we employ to approximate $\|L^\dagger\|_\infty$ (for Cases 2-4) and recall Hager's algorithm for estimating the 1-norm of the inverse of a matrix (for Case 1).
Subsequently, we describe our implementation based upon a boundary value problem formulation for nonautonomous inertial manifold reduction. Further details are available in \cite{CSVV14}.
In addition, we describe a method for obtaining local solutions to the fundamental matrix solution that does not require explicit formulas for the differential
equation or the linear variational equation on the inertial manifold. 

\subsection{Conditioning of IVP}

We recall Hager's algorithm to determine lower bounds on the $L^1$ condition number of a square matrix (see \cite{Hager}, \cite{Higham}).
For our purposes here we wish to estimate $\|L^{-1}\|$ for classical forward error analysis (Case 1) in which $L$ is a square matrix. The pseudo code for the code we employ is
Algorithm \ref{Hageralg}.

\begin{algorithm}[ht]\label{Hageralg}
        \caption{Hager's algorithm to approximate $\|L^{-1}\|_\infty$ in Case 1.}
        \label{alg:Hager}
        \begin{algorithmic}[1]
                \REQUIRE $\{X_n\}_{n=0}^N$
                \ENSURE Estimate (lower bound) $\|L^{-1}\|_\infty$
                \STATE $i_1 = -1; c_1=0; b$ random vector;
                \WHILE{(1)}
                \STATE{Solve $L^Tx = b$; Set $x = {\tt sign}(x)$;}
                \STATE{Solve $L b = x$; Set $i_2 = \max_i |b_i|$;}
                \IF{$1\leq i_1$}
                \IF{$i_1\equiv i_2$ {\bf or} $c_2 \leq c_1$}
                \STATE{break;}
                \ENDIF
                \ENDIF
                \STATE{$i_1=i_2; c_1=c_2; b=0; b_{i_1}=1;$}
                \ENDWHILE
        \end{algorithmic}
\end{algorithm}

We note that for classical forward error analysis (Case 1) $L$ is block unit lower triangular,
so no factorization is necessary when solving the linear systems in Hager's algorithm, only
block backward and forward substituion. Typically 4 linear system solves, two with $L$ and two with
$L^T$, are required.

In Cases 2-4 we consider two options for approximating $\|L^\dagger\|_\infty$. The first option which is the most expensive 
essentially involves computation of $\|L^\dagger|_\infty$ by iteratively determining the rows of $L^\dagger$ by computing the
columns of $(L^\dagger)^T$. This is accomplished by solving linear systems from $(LL^T)\cdot (L^\dagger)^T = L$ using the triadiagonal matrix $LL^T$.
This is done in a factor/solve framework by first performing a block tridiagonal factorization analogous to Thomas' algorithm
and then solving with multiple right hand sides. 
The second option relies on using a time stepping technique such as an embedded Runge-Kutta pair that provides an estimate
of $(G(x)_n = x_{n+1}-\varphi(x_n,h_n)$, for example by employing the higher order method as an approximation of $\varphi(x_n,h_n)$ with the lower order 
method providing $x_{n+1}$. We then solve the linear system $L(\Delta x, \Delta h) = -G(x)$, which is equivalent to finding the first Newton step,
and then approximating $\|L^\dagger\|_\infty \approx \|L^\dagger G(x)\|_\infty/\|G(x)\|_\infty$. This only requires a single linear system solve
and provides a residual correction based upon the approximation of the local error that is employed for step-size selection and error control.

Finally we note here the relationship between Case 3 (forward error analysis with rescaling of time) and Case 4 (shadowing error analysis with rescaling of time).
If we let $L_3$ denote the operator $L$ in Case 3 and $L_4$ for Case 4, then $L_4L_4^T = L_3 L_3^T + U U^T$ where 
$U^T = (X_0^T 0 \cdots 0)$. Then by the Sherman-Morrison-Woodbury formula
\[
(L_3L_3^T)^{-1} = (L_4 L_4^T)^{-1} + (L_4 L_4^T)^{-1} U (I - U^T(L_4 L_4^T)^{-1}U)^{-1}U^T (L_4 L_4^T)^{-1}
\]
provided the low dimensional matrix, the capacitance matrix, $I - U^T(L_4 L_4^T)^{-1}U$ is invertible. Note that $U^T (L_4 L_4^T)^{-1}$ is
the first block row of $L_4^\dagger$. From this relationship we see that the difference between the conditioning
of Cases 3 and 4 depends on the behavior of the capacitance matrix. Note also that this provides a means of updating from Case 3 to Case 4
(or vice versa).

Case 3 can also be related to Case 1 and solved via a forward substitution scheme. Since $\Delta x_0 = 0$ we have for $n=0,1,...,N-1$,
\[
\Delta x_{n+1} - \theta f_n \Delta h_n = b_n - X_n \Delta x_n \equiv g_n.
\]
Then solving $L_3 (\Delta x, \Delta h) = b$ in a minimum norm least squares sense is equivalent to solving the underdetermined low dimensional linear systems at
each step in a minimum norm sense, for example using the pseudo inverse of the matrix $[I | -\theta f_n]$.

\subsection{Inertial Manifold Reduction}

To put the original initial value problem \eqref{ivplinvar} into the form  \eqref{imform} we must have access to the orthogonal transformation $Q(t)$.  To find $Q(t)$ such that $Q^T(t)X(t) = R(t)$ where $R(t)$ is of the form \eqref{Rmat} we express $Q(t)$ as a product of Householder matrices $Q^T(t) = Q_p(t)\cdot \hdots \cdot Q_1(t)$ and follow the continuous technique.  For the details, see \cite{DV1,DV2}.\\

Since $Q_i(t)$ is a Householder matrix we can write $Q_i(t) = \left[\begin{array}{cc} I_{i-1} & 0 \\0 & P_i(t) \end{array} \right]$ where $I_{i-1}$ is the $i-1$ dimensional identity matrix and $P_i(t) = I -2 v_i(t) v_i(t)^T$ where $v_i(t) \in \mathbb{R}^{d-i+1}$ with $\|v_i\|_2 =1$.  To determine $v_i$, let $X_i$ be the transformed matrix $Q_{i-1}\cdot \hdots \cdot Q_1 X$. Then notice that if $u_i = X_i e_i -\sigma_i \|X_i e_i\|_2 e_i$, then $v_i = u_i/\|u_i\|_2$ defines a Householder transformation $Q_i$ that diagonalizes the $i^{th}$ column of $X(t)$.  The value $\sigma_i$ is chosen to ensure numerical stability and the canonical choice is 
(see e.g. \cite{GVL96})
\begin{equation}\label{numstab}
\sigma_i = \left\{\begin{array}{cc} -1 & \text{ if } e_1^T x_i(t) \geq 0 \\
 1 & \text{ if } e_1^T x_i(t) < 0 \end{array} \right.
\end{equation}
where $x_i(t)$ is the $i^{th}$ column of $X_i(t)$.  To further reduce the number of equations we need to solve we can define $w_i = v_i/(e_1^T v_i)$ and then notice that we can recover $v_i$ from $w_i$ via the additional relationship $e_1^T v_i = -\sigma_i / \|w_i\|_2$.  We want to avoid computing the fundamental matrix solution $X(t)$ or any of its columns and want to work only with the matrix $L(t)$ and the $w_i$'s.  To obtain $R(t)$ from $L(t)$, let $L_i(t)$ be the matrix obtained after diagonalizing the $i^{th}$ column of $X$.  We find can inductively find $L_i$ by doing a sequence of $(L,Q_i)$ updates:
\begin{equation}\label{LQudpate}
L_i(t) = Q_i(t)L_{i-1}(t)-Q_i(t)\dot{Q}_i(t), \quad i=1,\hdots,i-1
\end{equation}
where $L_0$ is taken to be $L$.  We can express the $(L,Q_i)$ update in terms of only the $L_i$ and $w_i$ as
$$Q_i(t)L_{i-1}(t)-Q_i\dot{Q}_i = L_{i-1}-\frac{2}{w_i^T w_i}\left( w_i (w_i^T L_{i-1})+(L_{i-1} w_i) w_{i}^T \right) = 4\frac{w_i^T L_{i-1} w_{i}}{(w_i^T w_i)^2}w_i w_i^T - \frac{2}{w_i^T w_i}(w_i \dot{w}_i^T-\dot{w}_i w_i^T).$$
From the definition of $w_i$ we have $w_i = (1, \hat{w}_i)^T$ and we can derive the following differential equation satisfied by $\hat{w}_i$:
\begin{equation}\label{weqn}
\dfrac{d \hat{w}_i}{dt} = \left(L_{i-1}^{1,1}]\hat{w}_i+\hat{w}_i^T \hat{L}^{1}_{i-1} - 2\frac{w_i^T L w}{w_i^T w_i}\right)\hat{w}_i +\left(1-\frac{w_i^T w_i}{2} \right) \hat{L}_{i-1}^{1} + \hat{L}_{i-1} \hat{w}_{i-1} \equiv h_i(\hat{w},t)
\end{equation}
where $L_{i-1}^{j,k}$ is the $(j,k)$ entry of $L_{i-1}$, $\hat{L}_{i-1}$ is the submatrix $L_{i-1}^{2:n,2:d}$ of $L_{i-1}$ and $\hat{L}^{1}_{i-1}$ is the column vector $L_{i-1}^{2:n,1}$.
The condition \eqref{numstab} is expressed by having 
\begin{equation}\label{wreembedcond}
1-\hat{w}_i^T \hat{w}_i \geq 0
\end{equation}
be satisfied at all times.  If \eqref{wreembedcond} is not satisfied, then we redefine $\sigma_i$ accordingly and then reembed the new update $\hat{w}_i$ variables to be consistent with the new $\sigma_i$.  Thus, starting from the original equation $\dot{u}=L(t)u+N(u,t)$, to find a point on the inertial manifold $x_0:=x(t_0)=\Psi(t_0,y(t_0)$ using the BVP formulation \eqref{bvpform} we use the following boundary value problem.
\begin{equation}\label{bvpQ}
\left\{\begin{array}{c}
\dot{x} = A(t)x + F(t,x,y)\\
\dot{y} = B(t)y+G(t,x,y)\\
\hat{w}_i = h_i(t,\hat{w}_{i-1}), \quad i=1,\hdots,p\\
x(-\infty) = 0, y(t_0) = y_0, \hat{w}_i(-\infty) = \hat{w}_{i,-\infty}
\end{array}
\right.
\end{equation}
In addition to \eqref{bvpQ} we have the initial value problem which is useful in constructing an initial guess to the solution of the boundary value problem
\begin{equation}\label{Qivp}
\left\{\begin{array}{c}
\dot{x} = A(t)x + F(t,x,y)\\
\dot{y} = B(t)y+G(t,x,y)\\
\hat{w}_i = h_i(t,\hat{w}_{i-1}), \quad i=1,\hdots,p\\
x(t_0) = x_0, y(t_0) = y_0, \hat{w}_i(t_0)=\hat{w}_{i,0}
\end{array}
\right.
\end{equation}

To solve for the $\hat{w}_i$ equations using Matlab's IVP and BVP solvers, we must modify the solvers so that we can reembed the $\hat{w}_i$ variables when the $\sigma_i$'s are changed.  To do so, we modify the ode45 code so that we perform a reembedding on $\hat{w}_i(t_n)$ for $i=1,\hdots,p$ just before ode45 computes the candidate numerical approximation to $\hat{w}_i(t_{n+1})$ for $i=1,\hdots,p$.  Similarly we modify the bvp4c code so that the reembedding happens after the Newton iteration converges to the candidate numerical approximation of $\hat{w}_i(t_{n_+1})$ for $i=1,\hdots,p$.  This leads to the following algorithm to compute the value of $x_0=\Psi(t_0,y(t_0))$.

  \begin{algorithm}\label{IMalg}
  \caption{(Algorithm to approximate $\Psi(t_0,y(t_0))$)}
        \begin{algorithmic}[1]
                \REQUIRE $w_i(-\infty)$, $y(-\infty)$, $T_{\infty}$
                \ENSURE $\Psi(t_0,y(t_0))$
                \STATE{ Set $\infty = T_{\infty}$.}
                \STATE{(Construct Initial guess for BVP) Solve \eqref{Qivp} on $[-T_{\infty},t]$ using modified ode45 with the initial conditions $y(-\infty)$, $x(-\infty)=0$, and $w_{i}(-\infty)$ for $i=1,\hdots,p$. }
                \STATE{(Solve BVP) Use the computed initial guess solution to solve \eqref{bvpQ} on the interval $[-\infty,t]=[-T_{\infty},t]$ with boundary conditions $x(-\infty)=0$, $y(t)=y_0$, and $w_i(-\infty)$ using the modified bvp4c and let $\Psi(t,y)$ be the output of $x(t)$ from the the boundary value solver.  }
                \end{algorithmic}
\end{algorithm}

The solution of the boundary value problem \eqref{bvpQ} produces a single point $x(t)=\Psi(t,y(t))$ on the inertial manifold.  To continue this solution and compute a trajectory on the inertial manifold, is equivalent to solving the initial value problem $\dot{y}=B(t)y+G(t,\Psi(t,y),y)$, $y(t_0)=y_0$.  To solve this problem, we must be able to evaluate $\Psi(t,y(t))$ which requires the solution of a boundary value problem of the form \eqref{bvpQ}.  With this in mind let  $y_{n+1} = \mathcal{H}(\Delta t,t_n,\{y_n\}_{k=0}^{n})$ be an $s$-stage, $k$-step explicit numerical method that approximates the initial value problem  $\dot{y}=B(t)y+G(t,\psi(t,y),y)$, $y(t_0)=y_0$.  We use the following algorithm to compute an approximate trajectory on the inertial manifold.

  \begin{algorithm}\label{IMalgts}
    \caption{(Algorithm to approximate the solution of $\dot y = B(t)y + G(t,\Psi(t,y),y)$)}
        \begin{algorithmic}[1]
                \REQUIRE $T,\Delta t$, $T_{\infty}$, $t_0$, $\hat{w}_{i,-T_{\infty}}$
                \ENSURE{ $\{y_n\}_{n=0}^{N}$, $\{t_n\}_{n=0}^{N}$}
                \STATE {Set $\infty=T_{\infty}$, $N=0$}
                 \WHILE{ $t_N < T$} 
                 \STATE{ Determine the step-size $(\Delta t)_N$}
                 \STATE{Set $t_{N+1} = t_{N}+(\Delta t)_N$}
                 \STATE{(Construct Initial guess for BVP) Solve \eqref{Qivp} on $[-\infty,t_N]=[t_N-T_{\infty},t_N]$ using modified ode45 with the initial conditions $y(t_N-T_{\infty})$, $x(t_N-T_{\infty})=0$, and $w_{i}(t_N-T_{\infty})\}$ for $i=1,\hdots,p$.}
                \STATE{(Solve BVP) Solve \eqref{bvpQ} on $[-\infty,t_N]=[t_N-T_{\infty},t_N]$ with the boundary conditions $x(t_N-T_{\infty})=0$, $y(t_N-T_{\infty})=y_N$, and $\hat{w}_i(t_N-T_{\infty})=\hat{w}_{i,-T_{\infty}}$ using modified bvp4c}
                \STATE{Set $y_{N+1} = \varphi(\Delta t,t_N,\{y_N\}_{k=0}^{n})$, $N=N+1$ (Repeating the above step to approximate stage values needed to evaluate the method)}
                \ENDWHILE
                
        \end{algorithmic}
\end{algorithm}

A drawback of using Algorithms 2 and 3 for computations on the inertial manifold is that it requires us to specify quantities that are not given as the initial data for an initial value problem.  In the standard initial value problem framework we are given the value of $u(0)=u_0$ which is the initial condition for the equation $\dot u = f(t,u)$.  To use Algorithms 2 and 3 we must specify the value of $p=\text{dim}(Y)$, $T_{\infty} \approx \infty$, $\hat{w}_i(-\infty)$, and $y(0)$.  Since $u(0)=Q(0)(x(0),y(0))^T$, we can only determine the value of $y(0)$ if we know the value of $Q(0)$ which requires knowledge of the values of $\hat{w}_i(0)$ which are only found from the solution of the boundary value problem in 2.  For each pair of $T_{\infty}\approx \infty$ and $p=\text{dim}(Y)$ there corresponds a different approximation to $u(t)=Q(t)(x(t),y(t)^T$.  This fact is explored more in the numerical results of Section 5.2.\\

We can use the output of Algorithm 3 to form $X_n$ and $f_n$ as follows.  Let $\{y_n\}_{n=0}^{N}$, $\{t_n\}_{n=0}^{N}$ be the output of Algorithm 3 and fix $\Delta x > 0$ and $\Delta t > 0$ small.  As noted in Section 3.1, for $i=1,\hdots,p$ we have the approximation
\begin{align*}
\frac{\partial}{\partial y} G_1(t,y)e_j & \approx \frac{1}{\Delta x}(G_1(t, y+(\Delta x) e_j)-G_1(t,y))\\ 
& \approx \frac{1}{(\Delta x) (\Delta t)}(\varphi(\tau;y+(\Delta x) e_j)-\varphi(0;y+(\Delta x) e_j)-\varphi(\Delta t;y)+\varphi(0;y))\\
& =\frac{1}{(\Delta x) (\Delta t)}\left(\varphi(\tau;y+(\Delta x) e_j)-\varphi(\Delta t;y)+(\Delta x) e_j \right)
\end{align*}
We first use the approximation $y_{n+1} \approx \varphi(\tau,y_n)$.  For $j=1,\hdots,p$, we can find quantities $(\Delta y)_{n,j} \approx  \varphi((\Delta t)_n,y_n+(\Delta x) e_j)$ using the Algorithm 2 with $t=t_n$ and $y(t_n)=y_n+(\Delta x) e_j$ and we find $x_n$ as the output of algorithm 2 with $t=t_n$ and $y(t_n)=y_n$.  Set $\delta_{n,j} = \frac{1}{(\Delta x)(\Delta t)}((\Delta y)_{n,j}-y_{n+1}+(\Delta x) e_j)$ for $j=1,\hdots,p$.  We then form an approximation to $X_n$ and $f_n$ as $X_n \approx I+(\Delta t)_n \Delta_n$ where $\Delta_n = [\delta_{n,1}|\hdots|\delta_{n,p}]$ and $f_n \approx G_1(t_n,x_n,y_n)$.

%

\section{Numerical Results}\label{Results}

In this section we present numerical results that show the 
amplification factors given by $\|L^\dagger\|$ for the four cases global error types considered here
for two example problems. The first is the classical Lorenz equations, a three dimensional system
of ODEs. The second is the Kuramoto-Sivashinsky equation which is known to have a low dimensional
inertial manifold.

\subsection{Lorenz '63}

Our first example is the classical Lorenz equation
$$\begin{pmatrix}\dot x\cr \dot y\cr \dot z\cr\end{pmatrix} = \begin{pmatrix}\sigma(y-x)\cr
                                                      \rho x -xz-y\cr
                                                      xy - \beta z\cr\end{pmatrix}
.$$
We consider the parameter values $\sigma = 10$, $\beta = 8/3$ and
$\rho=28$ and the initial condition $(x(0),y(0),z(0))=(0,1,0)$.

We employ Runge-Kutta $(4,5)$ to solve the
the nonlinear differential equation with local error control and adaptive step-size selection. 
Using the time mesh determined in this way the local solutions of the linear variational equation
(the $X_n$) are determined using the forward Euler method.


In Table \ref{newtable1} and Table \ref{newtable2} we report on values of $\|L^\dagger\|_\infty$ for the four different cases
for final times $T=1,10,100$ and absolute local error tolerances ${\tt tol}=10^{-4},10^{-6},10^{-8}$. 
As a check we also compared the values obtained by forming $L^\dagger$ using
block tridiagonal linear system solves with the matlab command {\tt norm(pinv(L),Inf)} and found
agreement in all cases to high precision. 

\begin{figure*}[ht]
\centering
\includegraphics[width=1.0\textwidth]{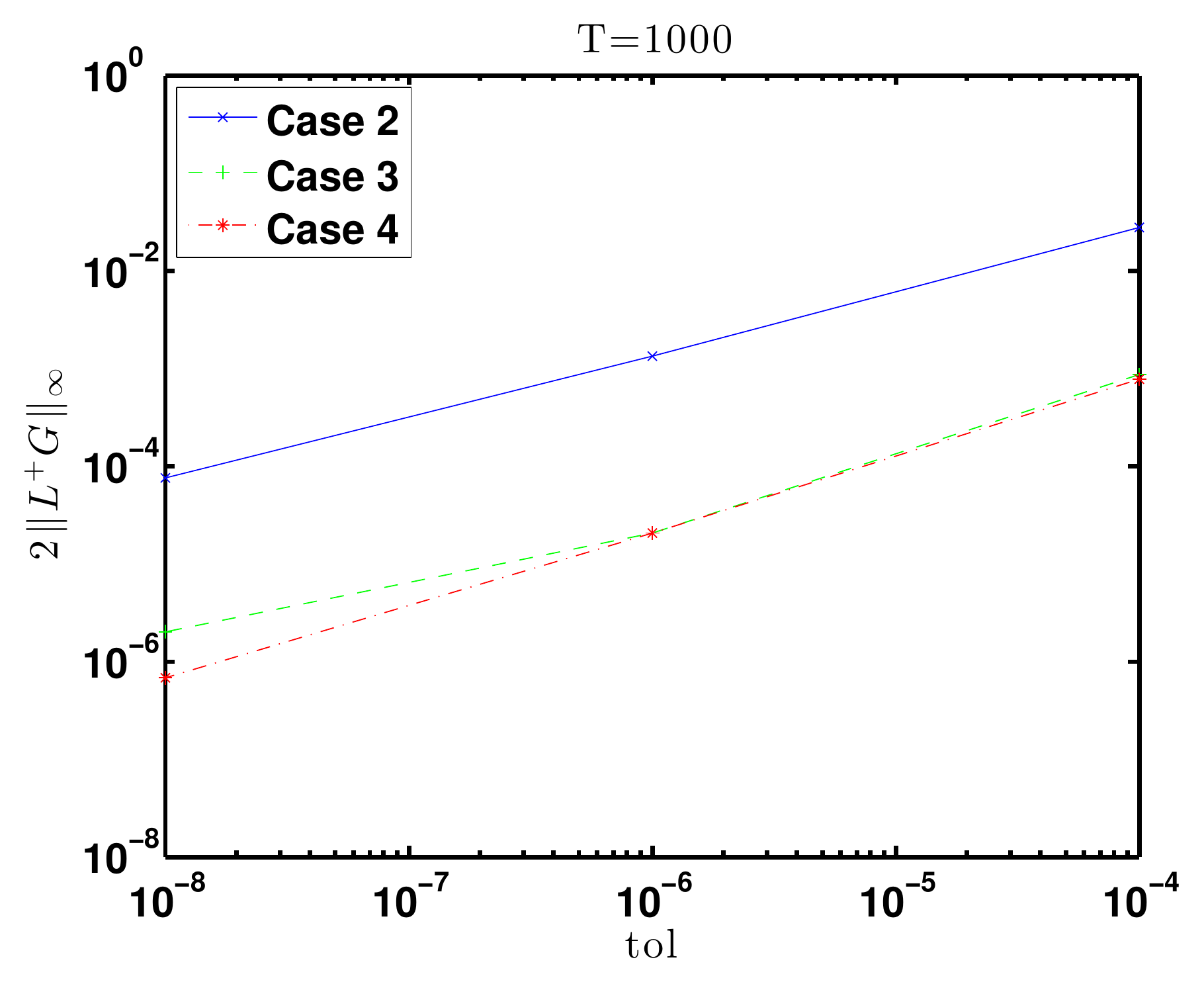}
\caption{Plot of ${\tt tol}$ vs. $2\|L^\dagger G\|$ for Cases 2,3,4 with $T=1000$.}
\label{figure1}
\end{figure*}
\begin{figure*}[ht]
\centering
\includegraphics[width=1.0\textwidth]{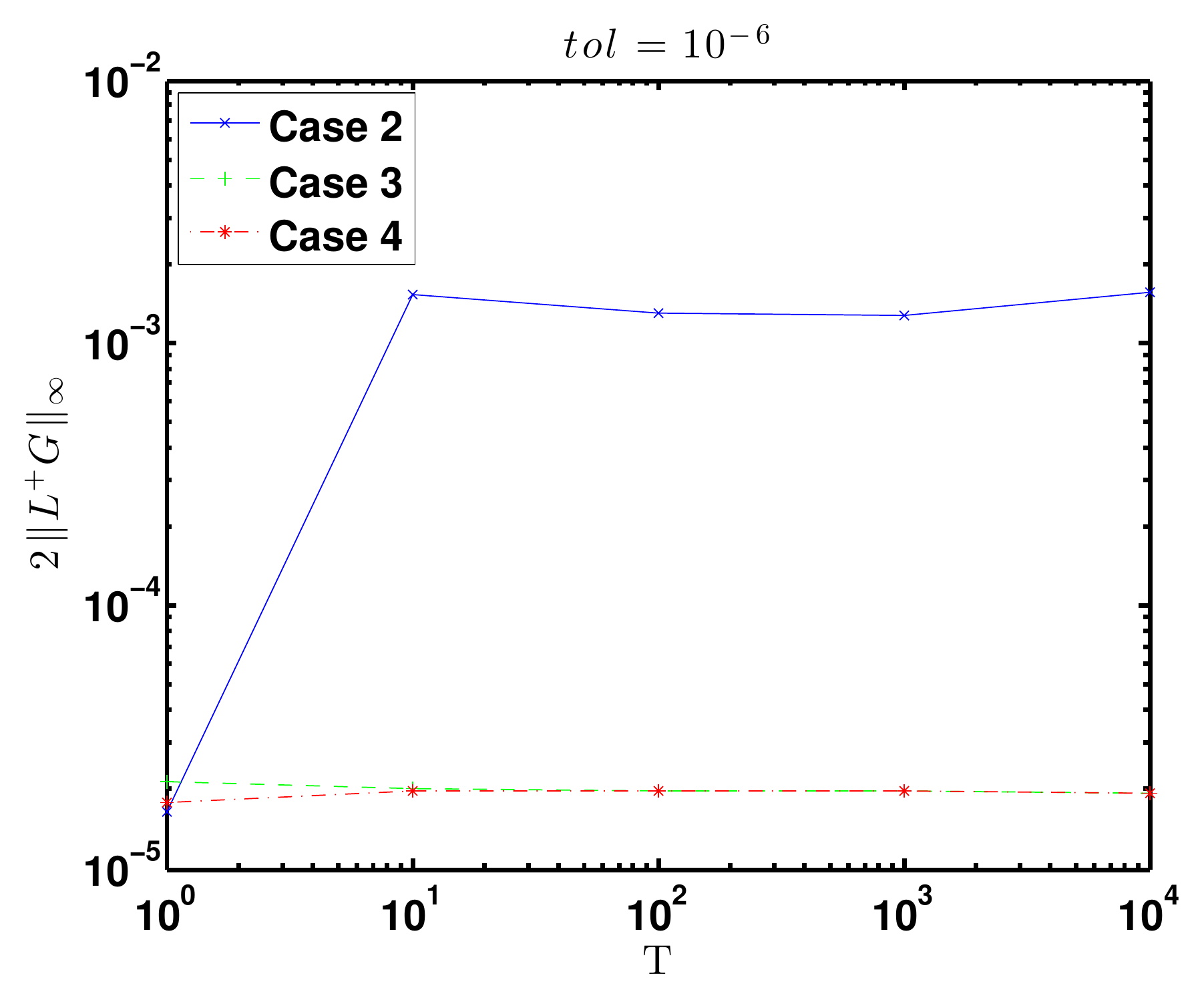}
\caption{Plot of $T$ vs. $2\|L^\dagger G\|$ for Cases 2,3,4 with ${\tt tol}=10^{-6}$.}
\label{figure2}
\end{figure*}

The value of $\theta$ used in Cases 3 and 4 controls the degree to which there is rescaling of time.
In general, the smallest possible value of $\theta>0$ is desired since the difference between the computed
time step and the time step of the exact solution, that is close to the computed solution, is bounded by
$2 \theta \|L^\dagger\|_\infty\cdot {\tt tol}$ when the nonlinear analysis holds. In the limit as $\theta\to 0$
Cases 3 and 4 revert to Cases 1 and 2, respectively, so we expect $\|L^\dagger\|_\infty$ to increase as $\theta\to 0$.
We determined ``optimal'' values of $\theta$ using final time $T=10^3$, ${\tt tol}=10^{-4}$ 
and found $\theta_{opt} \approx 60$ for Case 3 and $\theta_{opt} \approx 0.05$ for Case 4. We obtained similar
values of $\|L^\dagger\|$ for larger values of $\theta$ so these optimal values provide the tightest bound on the
time steps without significantly increasing $\|L^{\dagger}\|$. The value of the $\theta$ found in Case 4 agrees
with the value obtained in \cite{VVshad95} using a different time stepping technique. These are the values of
$\theta$ used to obtain the results in Table \ref{newtable2} and in Figures \ref{figure1} and \ref{figure2}.

In Figure \ref{figure1} we plot the tolerance ${\tt tol}$ versus the computed value of $\|L^\dagger G\|_\infty$ for Cases 2-4
for $T=1000$.
In Figure \ref{figure2} we plot the final time $T$ the computed value of $\|L^\dagger G\|_\infty$ for Cases 2-4
for ${\tt tol}=10^{-6}$. These plots reveal the good behavior of Case 3 (nearly as good as Case 4) and justify making
good long time global error claims in the sense of forward error analysis with rescaling of time.


\begin{table}
\begin{center}
\begin{tabular}{|c|c|c|c|c|c|}\hline
${\tt icase}$ & ${\tt tol}$ & $T$ & $\|L^{\dagger}\|_\infty$ & $\|L^{\dagger}G\|_\infty/\|G\|_\infty$ & {\tt Hager}\\ \hline\hline
 1 & $1E-4$ & 1 & $2.7E2$ & $1.6E2$ & $2.7E2$ \\ \hline
 1 & $1E-6$ & 1 & $5.5E2$ & $5.2E2$ & $5.5E2$ \\ \hline
 1 & $1E-8$ & 1 & $1.3E3$ & $1.3E3$ & $1.3E3$ \\ \hline\hline
 1 & $1E-4$ & 10 & $2.4E11$ & $7.5E10$ & $2.4E11$ \\ \hline
 1 & $1E-6$ & 10 & $5.3E6$ & $1.2E6$ & $5.3E6$ \\ \hline
 1 & $1E-8$ & 10 & $1.4E5$ & $3.8E4$ & $1.4E5$ \\ \hline\hline
 1 & $1E-4$ & 100 & $1.2E76$ & $1.0E75$ & $1.2E76$ \\ \hline
 1 & $1E-6$ & 100 & $2.1E50$ & $5.5E48$ & $2.1E50$ \\ \hline
 1 & $1E-8$ & 100 & $2.5E44$ & $4.1E42$ & $2.5E44$ \\ \hline\hline\hline
 2 & $1E-4$ & 1 & $5.0E1$ & $5.2E0$ & $--$ \\ \hline
 2 & $1E-6$ & 1 & $9.3E1$ & $8.2E0$ & $--$ \\ \hline
 2 & $1E-8$ & 1 & $2.1E2$ & $2.2E1$ & $--$ \\ \hline\hline
 2 & $1E-4$ & 10 & $2.1E3$ & $1.4E2$ & $--$ \\ \hline
 2 & $1E-6$ & 10 & $4.7E3$ & $7.7E2$ & $--$ \\ \hline
 2 & $1E-8$ & 10 & $1.2E4$ & $3.5E3$ & $--$ \\ \hline\hline
 2 & $1E-4$ & 100 & $2.1E3$ & $1.2E2$ & $--$ \\ \hline
 2 & $1E-6$ & 100 & $1.0E4$ & $6.4E2$ & $--$ \\ \hline
 2 & $1E-8$ & 100 & $1.3E4$ & $3.8E3$ & $--$ \\ \hline
\end{tabular}
\caption{Comparison of values of $\|L^\dagger\|_\infty$ with approximations of $\|L^\dagger\|_\infty$ for Cases 1 and 2.}\label{newtable1}
\end{center}
\end{table}

\begin{table}
\begin{center}
\begin{tabular}{|c|c|c|c|c|}\hline
${\tt icase}$ & ${\tt tol}$ & $T$ & $\|L^{\dagger}\|_\infty$ & $\|L^{\dagger}G\|_\infty/\|G\|_\infty$ \\ \hline\hline
 3 & $1E-4$ & 1 & $5.1E0$ & $3.0E1$ \\ \hline
 3 & $1E-6$ & 1 & $1.1E1$ & $9.2E1$ \\ \hline
 3 & $1E-8$ & 1 & $1.8E1$ & $2.6E2$ \\ \hline\hline
 3 & $1E-4$ & 10 & $3.7E1$ & $4.7E0$ \\ \hline
 3 & $1E-6$ & 10 & $3.0E2$ & $1.0E1$ \\ \hline
 3 & $1E-8$ & 10 & $2.9E3$ & $1.0E2$ \\ \hline\hline
 3 & $1E-4$ & 100 & $3.8E1$ & $4.3E0$ \\ \hline
 3 & $1E-6$ & 100 & $3.0E2$ & $1.0E0$ \\ \hline
 3 & $1E-8$ & 100 & $3.3E3$ & $1.0E2$ \\ \hline\hline\hline
 4 & $1E-4$ & 1 & $2.1E1$ & $3.3E0$ \\ \hline
 4 & $1E-6$ & 1 & $8.4E1$ & $8.9E0$ \\ \hline
 4 & $1E-8$ & 1 & $2.1E2$ & $1.5E1$ \\ \hline\hline
 4 & $1E-4$ & 10 & $2.5E1$ & $3.4E0$ \\ \hline
 4 & $1E-6$ & 10 & $1.8E2$ & $1.0E1$ \\ \hline
 4 & $1E-8$ & 10 & $9.5E2$ & $3.4E1$ \\ \hline\hline
 4 & $1E-4$ & 100 & $2.5E1$ & $3.8E0$ \\ \hline
 4 & $1E-6$ & 100 & $1.8E2$ & $9.9E0$ \\ \hline
 4 & $1E-8$ & 100 & $9.6E2$ & $3.4E1$ \\ \hline
\end{tabular}
\caption{Comparison of values of $\|L^\dagger\|_\infty$ with approximation of $\|L^\dagger\|_\infty$ for Cases 3 and 4.}\label{newtable2}
\end{center}
\end{table}


\subsection{KSE Equation}

As a second example consider the Kuramoto-Sivashinsky equation
in the form
\begin{equation}\label{eq:ourKSE}
\pd{\ut}{\tau} + 4\pd{^4\ut}{y^4} +
\vartheta\Bigl[\pd{^2\ut}{y^2} + \ut\pd{\ut}{y} \Bigr]= 0\;,
\end{equation}
with $\ut(y,t)=\ut(y+2\pi,t)$, and $\ut(y,t)=-\ut(-y,t)$.
With the change of variables
$$
-2w(s,y)=\ut(\xi s / 4,y), \qquad  \xi=\frac{4}{\vartheta}\;.
$$
\eqref{eq:ourKSE} can be written as
\begin{equation}\label{eq:CCP}
w_s=(w^2)_y-w_{yy}-\xi w_{yyyy}\;.
\end{equation}
All computations were performed on a spatial discretization using a standard Galerkin truncation with 16 modes (see \cite{DJRV} for more details) with the parameter values $\xi=0.02991$,  $\vartheta=133.73454$, which is one of the parameter values considered in \cite{DJRV,DJVV11}.  We present results in Table 4 for an approximate trajectory on the inertial manifold of \eqref{eq:ourKSE} found using Algorithm 3 and as a comparison in Table 3 we present results where the system \eqref{eq:ourKSE} is solved directly as in the case of the Lorenz '63 system.

For the direct solution of \eqref{eq:ourKSE} we employ a Runge-Kutta $(4,5)$ method with local error control and adaptive step-size selection.  For computations on the inertial manifold we employ the fixed time-step two-step Adams-Bashforth formula for the time-discretization of \eqref{eq:ourKSE} with time-step size given by $\Delta t=1E-3$.  To get an initial condition $(x(0),y(0))$ for the inertial manifold time-stepping scheme Algorithm 3 we fix $T_{\infty}=0.0005$, $p=\text{dim}(Y)$ and let the $j^{th}$ component of $y(0)$ be given by $(y_0)_j = (-1)^j/\sqrt{p}$ and $\hat{w}_{i,-\infty} = 0$ for $i=1,\hdots,p$. After this, Algorithm 2 is used to find the values of $x(0)$ and $\hat{w}_{i}(0)$ for $i=1,\hdots,p$ and $Q(0)$ is formed from the values of the $\hat{w}_{i}(0)$ to obtain $u(0)=Q(0)(x(0),y(0))^T$.  The local solutions of the linear variational equation (the $X_n$) are determined using the forward Euler method for both the direct and inertial manifold computations.  In Table 3 we present results for the value of $\|L^{\dagger}\|_{\infty}$ for the direct solution of \eqref{eq:ourKSE}.  In Table 4 we present results for the value of $\|L^{\dagger}\|_{\infty}$ for the flow computed from Algorithm 3 for $p=8$, $T_{\infty}=0.0005$, and $\Delta t = 1E-3$.  For both experiments we use final times of $T=0.1$ and $T=1.0$.

It is evident from Tables 3 and 4 that the global error as measured by $\|L^{\dagger}\|_{\infty}$ for the computations on the inertial manifold is much less than for the direct simulation in Cases 1-4.  This is at least partially explained by the way we compute the inertial manifold reduction.  Algorithm 3 allows us to reduce the dimension of \eqref{eq:ourKSE} from $16$ modes to $p$ modes by running time-stepping on an equation of the form $\dot y = B(t)y + G(t,\Psi(t,y),y)$ where we form an approximation $\tilde{\Psi}(t,y)$ to the inertial manifold equations $x=\Psi(t,y)$ using the Algorithm \eqref{IMalg}.  In essence, our inertial manifold time-stepping algorithm allows us to ignore the stiffest $16-p$ modes and avoid the large errors associated with the stiffest components.  The disadvantage of our time-stepping algorithm is that at each time-step we must solve a boundary value problem to compute $\tilde{\Psi}(t,y)$.  However, this disadvantage may be partially offset since we may be able to use a larger step-size since the reduced problem is less stiff than the full problem.

\begin{table}
\begin{center}
\begin{tabular}{|c|c|c|c|c|c|}\hline
${\tt icase}$ &  $T$ &  $\|L^{\dagger}\|_\infty$  \\ \hline\hline
      1      &  $0.1$ &   $5.0E79$                      \\\hline 
      1      &  $1$ &     $\text{Inf}$                    \\\hline \hline
      2      &   $0.1$  &  $2.4E1$                      \\\hline
      2      &   $1$  &    $2.5E1$                   \\\hline  \hline
      3      &  $0.1$  &  $3.2E21$              \\\hline 
      3      &  $1$ &     $3.2E33$                \\\hline \hline
      4      &   $0.1$  &  $2.0E1$                        \\\hline
      4      &    $1$  &    $2.5E1$                       \\\hline
\end{tabular}
\caption{Comparison of values of $\|L^\dagger\|_\infty$ for direct solution of a $16$ dimensional Galerkin approximation to \eqref{eq:ourKSE} for Cases 1-4.  The local error tolerance used was $1E-4$.}\label{KSEtable1}
\end{center}
\end{table}

\begin{table}
\begin{center}
\begin{tabular}{|c|c|c|c|c|c|}\hline
${\tt icase}$ &  $T$ &  $\|L^{\dagger}\|_\infty$  \\ \hline\hline
      1      &  $0.1$ &   $8.4E1$                      \\\hline 
      1      &  $1$ &     $1.0E3$                    \\\hline \hline
      2      &   $0.1$  &  $4.2E1$                      \\\hline
      2      &   $1$  &    $5.3E2$                   \\\hline  \hline
      3      &  $0.1$  &  $1.3E2$              \\\hline 
      3      &  $1$ &     $9.9E2$                \\\hline \hline
      4      &   $0.1$  &  $7.1E1$                        \\\hline
      4      &    $1$  &    $7.4E2$                       \\\hline
\end{tabular}
\caption{Comparison of values of $\|L^\dagger\|_\infty$ for Cases 1-4 for the approximate flow on the inertial manifold of the $16$ dimensional Galerkin approximation of \eqref{eq:ourKSE} using $\text{dim}(Y)=p=8$, time step-size $\Delta t = 1E-3$ and $T_{\infty} = 0.0005$.}\label{KSEtable2}
\end{center}
\end{table}

\section{Discussion}\label{Discuss}

In this paper we consider four types of global error analysis for initial value differential equations.
These can be thought of as different types of condition numbers for initial value differential equations
and correspond to the magnification factor of the local error that determines the global error. While computationally expensive it is shown that these condition numbers can be approximated in ways that give the correct general behavior that could be used as a measure of confidence or lack of confidence of a numerical simulation. The techniques are applied to the classical Lorenz '63 model and are exercised for a recently developed inertial manifold
dimension reduction method applied the the Kuramoto-Sivashinsky equation. Interestingly, forward error analysis with rescaling of time (Case 3) is quite effective for the Lorenz model but not for the (full) KSE model.  This appears to be stabilization that occurs due to the low dimension of the Lorenz model with a single positive Lyapunov exponent.  Additionally, the size of the condition numbers for the inertial form reduction of the KSE model seem to grow at a much slower rate than for the full unreduced model.

\bibliographystyle{plain}
\bibliography{decoup}

\end{document}